# Is there a "loophole" in Gödel's interpretation of his formal reasoning and its consequences?

Bhupinder Singh Anand

We formally define a "mathematical object" and "set". We then argue that expressions such as "(A$x$)$F(x)$", and "(E$x$)$F(x)$", in an interpretation M of a formal theory P, may be taken to mean "$F(x)$ is true for all $x$ in M", and "$F(x)$ is true for some $x$ in M", respectively, if, and only if, the predicate letter "$F$" is a mathematical object in P.

In the absence of such a meta-proof, the expressions "(A$x$)$F(x)$", and "(E$x$)$F(x)$", can only be taken to mean that "$F(x)$ is true for any given $x$ in M", and "It is not true that $F(x)$ is false for any given $x$ in M", respectively, indicating that the predicate "$F(x)$" is well-defined, and effectively decidable individually, for any given value of $x$, but that there may be no uniform effective method (algorithm) for such decidability.

We show how some paradoxical concepts of Quantum mechanics can then be expressed in a constructive interpretation of standard Peano's Arithmetic.

**1. Introduction**

Is there a "loophole" in Gödel's conclusion, in his seminal 1931 paper [Go31a], that his formally undecidable proposition GUS - which we may formally write as the string[1] [(A*x*)*R*(*x*)] - is true, under the standard interpretation M of Dedekind's standard Peano Arithmetic PA, in a constructive, and intuitionistically unobjectionable, way[2]?

---

[1] We use the synonymous terms "string" and "formula" interchangeably.

Now, classically, the string [(A*x*)*R*(*x*)] does translate formally as a true proposition under the standard interpretation of PA. However, there is nothing intuitive or constructive (in the sense of being effectively verifiable) about such "truth". The proposition is "true" only in a formal sense, where we accept Tarski's definition[3]: a string [(A*x*)*F*(*x*)] is true under an interpretation M if, and only if, it is satisfied by every *x* in M.

Although this appears to be a fairly innocent definition of intuitive truth, we note that it (implicitly) implies that we may (explicitly) assert a proposition as true in M if, and only if, it is severally, *and* jointly, satisfied by the elements in the ontology of M (even if there are elements of this ontology that are not interpretations of any mathematical objects in PA).

Now, a constructive view of Gödel's reasoning is that it actually establishes the invalidity of such a broad, implicit (hence ambiguous) assumption; thus, we can argue that we *cannot* unrestrictedly assert [*R*(*x*)] as jointly satisfied by all *x* in the standard interpretation M of PA, although we *can* assert that [*R*(*x*)] is satisfied (severally) by any given *x* in M.

The significance of this distinction is seen in Corollary 1.1 to Meta-lemma 1 of Anand [An02b]; the totality of values for which [*R*(*x*)] is satisfied in M may not be constructively definable as a formal mathematical object. In other words, the classical assumption that such values define a "set" in an Axiomatic Set Theory such as ZFC may introduce an inconsistency. Thus, we cannot make any constructive, or meaningful, assertion about the totality of values that satisfy [$R(x)$] in M.

---

[2] At the end of his proof of Theorem VI in his 1931 paper [Go31a] on formally undecidable propositions, Gödel remarks: "One can easily convince oneself that the proof we have just given is constructive (for all the existential assertions occurring in the proof rest upon Theorem V which, as it is easy to see, is intuitionistically unobjectionable), ...".

[3] We take Mendelson ([Me64], p49-52) as a standard exposition of Tarski's definitions of the "satisfiability" and "truth" of formal strings under a given interpretation; and of the basic foundational concepts of classical mathematical theory in general.

More precisely, Corollary 1.1 proves that every recursive number-theoretic relation does not well-define a (recursively enumerable[4]) sub-set of the natural numbers consistently in any Axiomatic Set Theory[5] that is a model for PA.

The above becomes clearer if we define a mathematical object and a set precisely as follows:

> **Definition** (*i*): A *primitive mathematical object* is any symbol for an individual constant, predicate letter, or a function letter, which is defined as a primitive symbol of a formal mathematical language.

> **Definition** (*ii*): A *formal mathematical object* is any symbol for an individual constant, predicate letter, or a function letter that is either a primitive mathematical object, or that can be introduced through definition into a formal mathematical language without inviting inconsistency.

---

[4] A recursively enumerable set is classically defined ([Me64], p250) as the range of some recursive number-theoretic function, and is implicitly assumed to be consistent with any Axiomatic Set Theory that is a model for P.

[5] Loosely speaking, we cannot give a set-theoretic definition, of Gödel's primitive recursive function $Sb(x\ v|Z(y))$, such that $\{x\ |\ x=Sb(y\ 19|Z(y))\}$ defines a set in any Axiomatic Set Theory with a Comprehension Axiom, or its equivalent, without introducing inconsistency.

Consequently, we cannot consistently assume that every recursive function formally defines a recursively enumerable set. It follows that we are unable to define a recursive set as a recursively enumerable set whose complement is also recursively enumerable. In some cases, there may be no such complement.

It further follows that, if $F(y)$ is an arithmetical function such that $F(k) = Sb(k\ 19|Z(k))$ for any given $k$, the assertion that the expression $\{x\ |\ x=F(y)\}$ defines a formal set by the Comprehension Axiom may require additional qualification.

> **Definition** (*iii*): A *mathematical object* is any symbol that is either a *primitive mathematical object*, or a *formal mathematical object*.
>
> **Definition** (*iv*): A *set* is the range of any function whose function letter is a mathematical object.

We can then argue, as a consequence of the theorem cited above, that expressions such as "$(Ax)F(x)$", and "$(Ex)F(x)$", in an interpretation M of a formal theory P, may be taken to mean[6] "$F(x)$ is true for all $x$ in M", and "$F(x)$ is true for some $x$ in M", respectively, if, and only if, the predicate letter "$F$" is a formal mathematical object in P.

In the absence of a proof that "$F$" is such a mathematical object, the expressions "$(Ax)F(x)$" and "$(Ex)F(x)$" can only be taken to mean that "$F(x)$ is true for any given $x$ in M", and "It is not true that $F(x)$ is false for any given $x$ in M", indicating that the predicate "$F(x)$" is well-defined, and decidable, for any given value of $x$, but that there may not be any uniformly effective method for such decidability.

## 2. We do not need to accept absolute limits on what we can formalise

What exactly does this mean, and why is the distinction important?

Taking the latter part of the question first, we note that the central issue in the development of AI is that of finding effective methods of duplicating the cognitive and expressive processes of the human mind. This issue is being increasingly brought into sharper focus by the rapid advances in the experimental, behavioral, and computer sciences[7]. Penrose's "The Emperor's New Mind" and "Shadows of the Mind" highlight what is striking about the attempts, and struggles, of current work in these areas to express their observations adequately - necessarily in a predictable way - within the standard interpretations of formal propositions as offered by classical theory.

---

[6] We note that such interpreted symbolic expressions are simply abbreviations of semantically well-defined assertions of a language of communication in which we express our cognitive experiences.

[7] See, for instance ([RH01], footnote 18).

So, the question arises: Are formal classical theories essentially unable to adequately express the extent and range of human cognition, or does the problem lie in the way formal theories are classically interpreted at the moment? The former addresses the question of whether there are absolute limits on our capacity to express human cognition unambiguously; the latter, whether there are only temporal limits - not necessarily absolute - to the capacity of classical interpretations to communicate unambiguously that which we intended to capture within our formal expression.

Now, my thesis is that we may comfortably reject the first, by recognising that we can, indeed, constructively reformulate Tarski's definitions, of the "satisfiability" and "truth" of formal propositions under a given interpretation, less ambiguously.

**3. A constructive definition of Tarskian truth**

So, we return to the first part of the earlier question: what exactly does this mean?

Well, it means that we can, for instance, replace the strong[8], classical, implicit interpretation of Tarski's non-constructive definition, which is essentially the assertion:

> (*i*) The string $[F(x)]$ of a formal system P is true under an interpretation M of P if, and only if, the interpreted predicate $F(x)$ is satisfied jointly, *and* severally, by the elements in the ontology of M;

by the weaker, explicit assertion:

---

[8] As remarked earlier, such a strong assertion may also be invalid under a constructive interpretation of Gödel's reasoning in Theorem VI of his 1931 paper.

(*ii*) The string [*F*(*x*)] of a formal system P is true under an interpretation M of P if, and only if, the interpreted predicate *F*(*x*) is satisfied *either* jointly, *or* severally, by the elements in the ontology of M;

We thus arrive at the constructive - in the sense of being effectively verifiable - definitions:

(*iii*) The string [*F*(*x*)] of a formal system P is *individually* true under an interpretation M of P if, and only if, given any value $k$ in M, there is an individually effective method (which may depend on the value $k$) to determine that the interpreted proposition *F*($k$) is satisfied in M;

(*iv*) The string [*F*(*x*)] of a formal system P is *uniformly* true under an interpretation M of P if, and only if, there is a uniformly effective method (necessarily independent of *x*) such that, given any value $k$ in M, it can determine that the interpreted proposition *F*($k$) is satisfied in M.

(*v*) The string [*F*(*x*)] of a formal system P is true under an interpretation M of P if, and only if, it is either uniformly true in M, or it is individually true in M.

Clearly, if [*F*(*x*)] is uniformly true in M, then it is, obviously, individually true in M. However, does the converse hold? More to the point, could Gödel's undecidable proposition, [(A*x*)*R*(*x*)], be an instance of a formula [*R*(*x*)] that is individually true in M (since Gödel shows that [*R*($k$)] is provable in P for any numeral [$k$]), but not uniformly true in M?

An interesting consequence of an affirmative answer to this question would be that the interpreted arithmetical predicate *R*(*x*) - which is instantiationally equivalent to a primitive recursive relation - becomes Turing-undecidable! Prima facie, this would appear to conflict with the classical postulation that a number-theoretic function is Turing-computable if, and only if, it is partial recursive[9]. However, the conflict may be illusory: the proof of equivalence seems to presume that there is a uniformly effective method (algorithm) for computing the partial recursive function[10]. This proof would be invalid if we held that Gödel has shown there are recursive functions, and relations, that may not be computable, or decidable, respectively, by any uniform method.

(This may, indeed, have been the interpretation that prompted Turing to assert the essential equivalence between Gödel's definition of an undecidable, but true, arithmetical proposition, and his own definition of a number-theoretic Turing-uncomputable Halting function.)

**4. Effective computability and Church's Thesis**

A key question, of course, is: why change?

One reason is that of verifiability: we are now able to define effective computability constructively, i.e. in a verifiable manner, by answering the questions:

>   (*vi*) When may we constructively assume that, given any sequence *s* of an interpretation M of P, there is an individually effective method to determine that *s* satisfies a given P-formula in M?

>   (*vii*) When may we constructively assume that there is a uniformly effective method such that, given any sequence *s* of an interpretation M, *s* satisfies a given P-formula in M?

If the domain D of M can be assumed representable in P, then (*vi*) can be answered constructively; we simply reformulate the classical Church Thesis as follows:

---

[9] Cf. ([Me64], p233, Corollary 5.13 and p237, Corollary 5.15).

[10] See ([Me64], p226, Corollary 5.11).

(*viii*) **Individual Church Thesis**: If, for a given relation $R(x)$, and any value $a$ in some interpretation M of P, there is an individually effective method such that it will determine whether $R(a)$ holds in M or not, then every element of the domain D of M is the interpretation of some term of P, and there is some P-formula $[R'(x)]$ such that:

$R(a)$ holds in M if, and only if, $[R'(a)]$ is P-provable.

In other words, the Individual Church Thesis postulates that the domain of a relation $R$ that is effectively decidable individually in an interpretation M of some formal system P can only consist of mathematical objects, even if $R$ is not, itself, a mathematical object.

However, (*vii*) can be answered constructively for any interpretation M of P, if we postulate:

(*ix*) **Uniform Church Thesis**: If, in some interpretation M of P, there is a uniformly effective method such that, for a given relation $R(x)$, and any value $a$ in M, it will determine whether $R(a)$ holds in M or not, then $R(x)$ is the interpretation in M of a P-formula $[R(x)]$, and:

$R(a)$ holds in M if, and only if, $[R(a)]$ is P-provable.

Thus, the Universal Church Thesis postulates, firstly, that the domain of a relation $R$ that is effectively decidable uniformly in an interpretation M of a formal system P can only consist of mathematical objects; and, secondly, that $R$, too, is necessarily a mathematical object.

Some interesting consequences of the above are that:

(*x*) The Uniform Church Thesis implies that a formula $[R]$ is P-provable if, and only if, $[R]$ is uniformly true in some interpretation M of P.

(*xi*) The Uniform Church Thesis implies that if a number-theoretic relation $R(x)$ is uniformly satisfied in some interpretation M of P, then the predicate letter "$R$" is a formal mathematical object in P (i.e. it can be introduced through definition into P without inviting inconsistency).

(*xii*) The Uniform Church Thesis implies that, if a P-formula [$R$] is uniformly true in some interpretation M of P, then [$R$] is uniformly true in every model of P.

(*xiii*) The Uniform Church Thesis implies that if a formula [$R$] is not P-provable, but [$R$] is classically true under the standard interpretation, then [$R$] is individually true, but not uniformly true, in the standard model of P.

(*xiv*) The Uniform Church Thesis implies that Gödel's undecidable sentence GUS is individually true, but not uniformly true, in the standard model of P.[11]

By defining effective computability, both individually and uniformly, along similar lines, we can give a constructive definition of uncomputable number-theoretic functions:

(*xv*) A number-theoretic function $F(x_1, ..., x_n)$ in the standard interpretation M of P is uncomputable if, and only if, it is effectively computable individually, but not effectively computable uniformly.

This, last, removes the mysticism behind the fact that we can define a number-theoretic Halting function that is, paradoxically, Turing-uncomputable.

## 5. The Uniform Church Thesis and the classical Church-Turing Theses

The significance, of defining the truth of a formula of P under interpretation explicitly in terms of individual, and uniform, effective methods, and of expressing Church's Thesis constructively, is seen if we note that any computer can be designed to recognise a "looping" situation; it simply records every instantaneous tape description at the execution of each machine instruction, and compares the current instantaneous tape description with the record.

---

[11] An intriguing consequence of this argument is considered in Appendix 1.

Now, we could instruct such a machine to assign arbitrary values to those undefined instances of the Halting predicate whose occurrences cause the machine to loop. Prima facie, what we appear to have here is an individually effective decision method that is dependent entirely on the particular Halting function that is being computed, and cannot be predicted.

Is such a machine a Turing machine? The classical answer to this question is not obvious if we do not appeal to the Church-Turing Thesis.

However:

>   (*xvi*) If we assume a Uniform Church Thesis, then every partial recursive number-theoretic function $F(x_1, ..., x_n)$ has a unique constructive extension as a total function.

>   (*xvii*) If we assume a Uniform Church Thesis, then not every effectively computable function is classically Turing computable (so Turing's Thesis does not, then, hold).

>   (*xviii*) If we assume a Uniform Church Thesis, then not every (partially) recursive function is classically Turing-computable.[12]

>   (*xix*) If we assume a Uniform Church Thesis, then the class P of polynomial-time languages in the P versus NP problem may not define a formal mathematical object.

---

[12] The classical proof that every (partially) recursive function is classically Turing-computable uses induction over (partial) recursive functions, thus assuming that every such function is a mathematical object; by Meta-lemma 1, such an assumption is invalid.

**Appendix 1: Constructivity and classical Quantum Theory**

In Anand [An02c], we show how the introduction of constructive definitions of classical mathematical concepts may permit formal systems of Peano Arithmetic to model some of the more paradoxical concepts of Quantum Mechanics. For instance, as a consequence of §4(*xiv*), consider the following argument:

(*a*) Gödel has proved in his 1931 paper that there is an arithmetic formula $[R(x)]$ such that, for any given $k$, $[R(k)]$ is provable.

(*b*) Hence, for any given $k$, there is always some effective method for evaluating the arithmetic expression $R(k)$.

(*c*) Gödel has also proved in the above paper that $[(Ax)R(x)]$ is not provable.

(*d*) *Thesis*: There is no uniform effective method (algorithm/Turing machine) that can evaluate the arithmetic expression $R(x)$ for any given $x$.

(*e*) Thus, $R(n)$ is individually computable, but not uniformly computable.

(*f*) *Theorem (provable by induction)*: For any given $k$, we can always find some effective method (algorithm/Turing machine) $T(k)$ that can compute $R(n)$ for all $n<k$, i.e. $T(k)$ terminates for all $n<k$, but it "loops" on input $k$. (Note: All methods that evaluate $R(n)$ for all $n<k$ cannot be non-terminating on input $k$; this would imply that $R(k)$ is undefined, which would contradict (*b*).)

(*g*) *Quantum interpretation*: The process of finding $T(k+1)$ can be corresponded, firstly, to the act of finding a suitable method of measuring the value $R(k)$ precisely, and, secondly, to the collapse of the wave function at $k$ as a result of the measurement; we then have the new "state" $T(k')$, which can evaluate the value of $R(n)$ for all $n<k'$, where $k<k'$, but not beyond!

(*h*) If, now, we have some law that determines the state $T(k')$ from the state $T(k)$ and the interaction at $k$, we have a deterministic interaction that is, nevertheless, absolutely unpredictable, where we may then define free will as absolute unpredictability. (We note that, if $k' > k+1$, we have a language that admits interactions that can leave the state $T(k')$ unchanged.)

Now we note that a counter-thesis to (*d*) would be:

(*i*) *Counter-Thesis*: There is some uniform effective method (algorithm/Turing machine) that can evaluate the arithmetic expression $R(x)$ for any given $x$.

It follows from Gödel's reasoning in [Go31a] that both (*d*) and (*i*) are effectively unverifiable, since they cannot be proved formally. We thus have two standard models of Peano Arithmetic - classical and constructive - that are mutually inconsistent. If we assume that both are consistent, the above argument indicates the interpretation that implies (*d*) may be the more suitable language for expressing concepts of classical Quantum Theory.

<*Web page*>: http://www.abelard.org/turpap2/tp2-ie.asp - index>

(*Acknowledgements: My thanks to Ms Shuchi Mehta for requesting a note on my perceived "loophole" in Gödel's 1931 paper. This paper is an edited version of that note.*)

(*Updated: Saturday 10th May 2003 7:18:00 AM IST by* re@alixcomsi.com)